\documentclass[12pt,reqno,twoside,a4paper]{article}
\usepackage {longtable}
\usepackage{amsmath,amssymb,amsfonts}
\usepackage[T2A]{fontenc}
\usepackage[cp866]{inputenc}
\usepackage[english]{babel}
\usepackage{amsthm}
\usepackage[mathscr]{eucal}
\def\{{\protect\lbrace}
\def\}{\protect\rbrace}

\newcommand{\Hom}{\operatorname{Hom}}

\newcommand{\Ker}{\operatorname{Ker}}

\newcommand{\Sing}{\operatorname{Sing}}

\begin{document}

\begin{center}
\textbf{Injective and Automorphism-Invariant Non-Singular Modules}
\end{center}
\begin{center}
A.A. Tuganbaev\footnote{National Research University "MPEI", Lomonosov Moscow State University,\\ e-mail: tuganbaev@gmail.com.}
\end{center}

\textbf{Abstract.} 
Every automorphism-invariant right non-singular $A$-module is injective if and only if the factor ring of the ring $A$ with respect to its right Goldie radical is a right strongly semiprime ring.

The study is supported by Russian Scientific Foundation (project 16-11-10013).

{\bf Key words:} automorphism-invariant module, injective module, quasi-injective module, strongly semiprime ring  

\textbf{1. Introduction and preliminaries}

All rings are assumed to be associative and with zero identity element; all modules are unitary.
A module $M$ is said to be \textit{automorphism-invariant} if $M$ is invariant under any automorphism of its injective hull. In \cite{DicF69}, Dickson and Fuller studied automorphism-invariant modules, when the underlying ring is a finite-dimensional algebra over a field with more than two elements. In \cite[Theorem 16]{ErSS13}, Er, Singh and Srivastava proved that a module $M$ is automorphism-invariant if and only if $M$ is a {\it
pseudo-injective} module, i.e.,  for any submodule $X$ of $M$, every monomorphism $X\to M$ can be extended to an endomorphism of the module $M$.
Pseudo-injective modules were studied in several papers; e.g., see \cite{JaiS75}, \cite{Tep75}, \cite{ErSS13}. Automorphism-invariant modules were studied in several papers; e.g., see \cite{AlaEJ05}, \cite{ErSS13}, \cite{GuiS13}, \cite{LeeZ13}, \cite{SinS13}, \cite{Tug13b}, \cite{Tug13d}, \cite{Tug14}, \cite{Tug15}, \cite{Tug17}. 

A module $X$ is said to be \textit{injective relative to the module} $Y$ or \textit{$Y$-injective} if for any submodule $Y_1$ of $Y$, every homomorphism $Y_1\to X$ can be extended to a homomorphism $Y\to X$. A module is said to be \textit{injective} if it is injective with respect to any module. A module $X$ is said to be {\it quasi-injective} if $X$ is injective relative to $X$. Every quasi-injective module is an automorphism-invariant module, since the module $X$ is quasi-injective if and only if $X$ is invariant under any endomorphism of its injective hull; e.g., see \cite[Theorem 6.74]{Lam99}. Every finite cyclic group is a quasi-injective non-injective module over the ring $\mathbb{Z}$ of integers.

A submodule $Y$ of the module $X$ is said to be \textit{essential} in $X$ if $Y\cap Z\ne 0$ for any non-zero submodule $Z$ of $X$. 
A submodule $Y$ of the module $X$ is said to be \textit{closed} in $X$ if $Y=Y'$ for every submodule $Y'$ of $X$ which is an essential extension of the module $Y$.
A module $M$ is called a \textit{CS module} if every its closed submodule is a direct summand of the module $M$. 

We denote by $\Sing X$ the \textit{singular submodule} of the right $A$-module $X$, i.e., $\Sing X$ is a fully invariant submodule of $X$ which consists of all elements $x\in X$ such that $r(x)$ is an essential right ideal of the ring $A$. A module $X$ is said to be \textit{non-singular} if $\Sing X=0$. 
For a module $X$, we denote by $G(X)$ the intersection of all submodules $Y$ of the module $X$ such that the factor module $X/Y$ is non-singular. The submodule $G(X)$ is a fully invariant submodule of $X$; it is called the {\it Goldie radical} of the module $X$. The relation $G(X)=0$ is equivalent to the property that the module $M$ is non-singular.

\textbf{Remark 1.1.} A ring $A$ is said to be \textit{right strongly semiprime} \cite{Han75} if any its ideal, which is an essential right ideal, contains a finite subset with zero right annihilator. The direct product of an infinite set of fields is an example of a commutative semiprime non-singular ring which is not strongly semiprime. All finite direct products of rings without zero-divisors and all finite direct products of simple rings are right and left strongly semiprime

Kutami and Oshiro \cite{KutO80} proved the following theorem.

\textbf{Theorem 1.2 \cite[Theorem 3.4]{KutO80}.} \textit{Every non-singular quasi-injective right module over ring $A$ is injective if and only if $A/G(A_A)$ is a right strongly semiprime ring.}

In connection to Theorem 1.2, we will prove Theorem 1.3 which is the first main result of the given paper.

\textbf{Theorem 1.3.} \textit{For a ring $A$ with right Goldy radical $G(A_A)$, the following conditions are equivalent.}
\begin{enumerate}
\item[\textbf{1)}] 
\textit{Every automorphism-invariant right non-singular $A$-module is an injective module.}
\item[\textbf{2)}] 
\textit{$A/G(A_A)$ is a right strongly semiprime ring.}
\end{enumerate}

\textbf{Remark 1.4.} Since every finite direct sum of injective modules and every direct product of injective modules are injective modules, it follows from Theorem 1.3 that the condition 2) of Theorem 1.3 is equivalent to each of the following conditions:\\
\textbf{3)} every finite direct sum of automorphism-invariant right non-singular $A$-modules is an injective module;\\
\textbf{4)} every direct product of automorphism-invariant right non-singular $A$-modules is an injective module.

The proof of Theorem 1.3 is decomposed into a series of assertions, some of which are of independent interest. 
We give some necessary definitions. 

In the paper, we use well-known properties of $\Sing X$, $G(X)$, and non-singular modules; e.g., see \cite[Chapter 2]{Goo76}, \cite[Section 7]{Lam99}, and \cite[Section 3.3]{Row88}. A module $Q$ is called an \textit{injective hull of the module} $M$ if $Q$ is an injective module and $M$ is an essential submodule of the module $Q$. 
A module is said to be \textit{square-free} if it does not contain the direct sum of two non-zero isomorphic submodules.
A module $M$ is said to be \textit{uniform} if the intersection of any two non-zero submodules of the module $M$ is not equal to zero. A module $M$ is said to be \textit{finite-dimensional} if $M$ does not contain an infinite direct sum of non-zero submodules.

\textbf{2. The proof of Theorem 1.3}

\textbf{Lemma 2.1.} \textit{Let $A$ be a semiprime ring containing a right ideal $D$ such that $D_A=\oplus _{i=1}^{\infty}D_i$, where each of the right ideals $D_i$ is non-zero and $D_iD_j=0$ for all $i\ne j$. Then the infinite sum of each of the two-sided ideals $AD_i$ of the ring $A$ is their direct sum.}

\textbf{Proof.} If $i\ne j$, then $(AD_i)(AD_j)=AD_iD_j=0$. In addition, 
$$
((AD_i)\cap (\sum_{j\ne i}(AD_j))^2\subseteq (AD_i)(\sum_{j\ne i}(AD_j))= (\sum_{j\ne i}(AD_i)(AD_j))=0.
$$
Since $A$ is a semiprime ring, $(AD_i)\cap (\sum_{j\ne i}(AD_j)=0$. Therefore, the infinite sum of each of the ideals $AD_i$ is their direct sum.\hfill$\square$

For convenience, we give the proof of the following familiar lemma.

\textbf{Lemma 2.2.} {\it Let $A$ be a ring and let $M$ be a non-zero right non-singular $A$-module. Then there exists a
non-zero right ideal $B$ of the ring $A$ such that the module $B_A$ is isomorphic to a submodule of the module $M$.}

\textbf{Proof.} Since $M$ is a non-zero non-singular module, there exists an element $m\in M$ such that $mA$ is a non-zero non-singular module. Since $mA\cong A_A/r(m)$ and the module $mA$ is non-singular, the right ideal $r(m)$ is not essential. Therefore, there exists a non-zero right ideal $B$ with $B\cap r(m)=0$. In addition, there exists an epimorphism $f\colon A_A\to mA$ with kernel $r(m)$. Since $B\cap \Ker f=0$, we have that $f$ induces the monomorpism $g\colon B\to mA$. Therefore, $bA$ contains the non-zero submodule $g(B)$ which is isomorphic to the module $B_A$.~\hfill$\square$

The following lemma is proved in \cite[Theorem 3, Theorem 6(ii)]{ErSS13}. 

\textbf{Lemma 2.3 \cite[Corollary 15]{LeeZ13}.} Let $M$ be the direct sum of CS modules $M_i$, $i\in I$. The module $M$ is quasi-injective if and only if $M$ is an automorphism-invariant module.

\textbf{Lemma 2.4.} \textit{Let $A$ be a ring and let $M$ be an automorphism-invariant right non-singular $A$-module.} 
\begin{enumerate}
\item[\textbf{1.}] 
\textit{$M = X\oplus Y$, where $X$ is a quasi-injective non-singular module, $Y$ is an automorphism-invariant non-singular square-free module, $\Hom (D_1,D_2) = 0$ for any two submodules $D_1$, $D_2$ of the module $Y$ with $D_1\cap D_2 = 0$. In addition, for any set $\{K_i\,|\, i\in I\}$ of closed submodules in $Y$, the submodule $\sum _{i\in I}K_i$ is automorphism-invariant.}
\item[\textbf{2.}] 
\textit{If $A$ is a semiprime ring which does not contain the infinite direct sum of non-zero two-sided ideals, then 
either $Y=0$ or the module $Y$ contains a non-zero quasi-injective uniform submodule.}
\end{enumerate}

\textbf{Proof.} \textbf{1.} The assertion is proved in \cite[Theorem 3, Theorem 6(ii)]{ErSS13}.

\textbf{2.} By Lemma 2.2, there exists a non-zero right ideal $B$ of the ring $A$ such that the module $B_A$ is isomorphic to a submodule of the module $Y$. 

We assume that $Y$ contains some non-zero uniform submodule. Then $Y$ contains some non-zero closed uniform submodule $Y_1$. By 1, $Y_1$ is an automorphism-invariant module. Since $Y_1$ is an automorphism-invariant uniform module and every uniform module is a CS module, it follows from Lemma 2.3 that $Y_1$ is a quasi-injective module, which is required.

Now we assume that $Y$ is a non-zero module which does not contain a non-zero uniform submodule. Then the right ideal $B$ contains the direct sum of $D_A=\oplus _{i=1}^{\infty}D_i$, where each of the right ideals $D_i$ is non-zero. For each $d_j\in D_j$, the mapping $f(x)=d_jx$ defines an $A$-module homomorphism $f\colon D_i\to D_j$. By 1, $\Hom (D_i,D_j) = 0$ for all $i\ne j$. Therefore, $D_iD_j=0$ for all $i\ne j$. By Lemma 2.1, the infinite sum of all the two-sided ideals $AD_i$ is their direct sum. By assumption, $A$ does not contain an infinite direct sum of non-zero two-sided ideals; this is a contradiction.\hfill$\square$

\textbf{Lemma 2.5.} \textit{Let $A$ be a semiprime ring which does not contain an infinite direct sum of non-zero two-sided ideals and let $M$ be an automorphism-invariant right non-singular $A$-module.} 
\begin{enumerate}
\item[\textbf{1.}] 
\textit{If every quasi-injective uniform submodule of the module $M$ is injective, then $M = X\oplus Y$, where $X$ is a quasi-injective non-singular module, $Y$ is an automorphism-invariant non-singular square-free module, and either $Y=0$ or $Y$ is an essential extension of some quasi-injective module $K$ which is direct sum of non-zero injective uniform modules.}
\item[\textbf{2.}] 
\textit{If every quasi-injective submodule of the module $M$ is injective, then $M$ is an injective module.}
\end{enumerate}

\textbf{Proof.} \textbf{1.} It follows from the assumption and Lemma 2.4 that $M = X\oplus Y$, where $X$ is a quasi-injective non-singular module, $Y$ is an automorphism-invariant non-singular square-free module, and either $Y=0$ or every non-zero direct summand of the module $Y$ contains a non-zero injective uniform direct summand. It is directly verified that $Y$ is an essential extension of some module $K$ which is the direct sum of non-zero injective uniform direct summands $K_i$, $i\in I$, of the automorphism-invariant non-singular square-free module $Y$. By Lemma 2.4(1), $K=\oplus_{i\in I}K_i$ is an automorphism-invariant module. By Lemma 2.3, $K$ is a quasi-injective module. 

\textbf{2.} For the module $M$, the assertion 1 of this lemma is true. In addition, every quasi-injective submodule of the module $M$ is injective. Therefore, the modules $X$ and $K$ from the condition 1 are injective, and $M=X\oplus Y$ and either $Y=0$ or $Y$ is an essential extension of the injective module $K$. If $Y=0$, then $M=X$ is an injective module. If $Y\ne 0$, then $Y=K$ and $M=X\oplus Y$  is an injective module.~\hfill$\square$

\textbf{Lemma 2.6 \cite{Han75}.} \textit{If $A$ is a right strongly semiprime ring, then $A$ is a right non-singular semiprime ring which does contain the infinite direct sum of non-zero ideals.} 

\textbf{Theorem 2.7.} \textit{For a ring $A$, the following conditions are equivalent.}
\begin{enumerate}
\item[\textbf{1)}] 
\textit{$A$ is a right non-singular ring and every automorphism-invariant right non-singular $A$-module is injective module.}
\item[\textbf{2)}] 
\textit{$A$ is a right strongly semiprime ring.}
\end{enumerate}

\textbf{Proof.} The implication 1)\,$\Rightarrow$\,2) follows from Theorem 1.2 and Lemma 2.5(2).

The implication 2)\,$\Rightarrow$\,1) follows from Theorem 1.2 and Lemma 2.6.~\hfill$\square$

\textbf{Lemma 2.8.} \textit{Let $A$ be a ring, $G=G(A_A)$ be the right Goldie radical of the ring $A$, $h\colon A\to A/G$ be the natural ring epimorphism, and let $X$ be a non-singular non-zero right $A$-module.}
\begin{enumerate}
\item[\textbf{1)}]
\textit{If $B$ is an essential right ideal of the ring $A$, then $h(B)$ is an essential right ideal of the ring $h(A)$.} 
\item[\textbf{2)}] 
\textit{If $B$ is a right ideal of the ring $A$ such that $G\subseteq B$ and $h(B)$ is an essential right ideal of the ring $h(A)$, then $B$ is an essential right ideal of the ring $A$.} 
\item[\textbf{3)}] 
\textit{For each right $A$-module $M$, the module $MG$ is contained in the Goldie radical of the module $M$.}
\item[\textbf{4)}] 
\textit{$XG=0$ and the natural $h(A)$-module $X$ is non-singular. In addition, if $Y$ is an arbitrary right non-singular $A$-module, then $YG=0$ and the $h(A)$-module homomorphisms $Y\to X$ coincide with the $A$-module homomorphisms $Y\to X$. Consequently, $X$ is an $Y$-injective $A$-module if and only if $X$ is an $Y$-injective $h(A)$-module. The essential submodules of the  $h(A)$-module $X$ coincide with essential submodules of the $A$-module $X$.}
\item[\textbf{5)}] 
\textit{$X$ is an injective $h(A)$-module if and only if $X$ is an injective $A$-module.}
\end{enumerate}

\textbf{Proof.} \textbf{1.} We assume that $h(B)$ is not an essential right ideal of the ring $h(A)$. Then there exists a right ideal $C$ of the ring $A$ such that $C$ properly contains $G$ and $h(B)\cap h(C)=h(0)$. Since $h(B)\cap h(C)=h(0)$, we have  $B\cap C\subseteq G$. Since $C$ strongly contains the closed right ideal $G$, we have that $C_A$ contains a non-zero submodule $D$ with $D\cap G=0$. Since $B$ is an essential right ideal, $B\cap D\ne 0$ and $(B\cap D)\cap G=0$. Then $h(0)\ne h(B\cap D)\subseteq h(B)\cap h(C)=h(0)$. This is a contradiction.

\textbf{2.} We assume that $B$ is not an essential right ideal of the ring $A$. Then $B\cap C=0$ for some non-zero right ideal $C$ of the ring $A$ and $G\cap C\subseteq B\cap C=0$. Therefore, $h(C)\ne h(0)$. Since $h(B)$ is an essential right ideal of the ring $h(A)$, we have $h(B)\cap h(C)\ne h(0)$. Let $h(0)\ne h(b)=h(c)\in h(B)\cap h(C)$, where $b\in B$ and $c\in C$. Then $c-b\in G\subseteq B$. Therefore, $c\in B\cap C=0$, whence we have $h(c)=h(0)$. This is a contradiction.

\textbf{3.} For each element $m\in M$, the module $mG_A$ is Goldie-radical, since $mG_A$ is a homomorphic image of the  Goldie-radical module $G$. Therefore, $mG\subseteq G(M)$ and $MG\subseteq G(M)$.

\textbf{4.} By 3, $XG=0$. We assume that $x\in X$ and $xh(B)=0$ for some essential right ideal $h(B)$, where $B=h^{-1}(h(B))$ is the complete pre-image of $h(B)$ in the ring $A$. By 2), $B$ is an essential right ideal of the ring $A$. Then $xB=0$ and $x\in \Sing X=0$. Therefore, $X$ is a non-singular $h(A)$-module.
The remaining part of 4 is directly verified.

\textbf{5.} Let $R$ be one of the rings $A$, $h(A)$ and let $M$ be a right $R$-module. It follows from the well-known Baer criterion that the module $M$ is injective if and only if $M$ is injective with respect to the module $R_R$. Now the assertion follows from 4.~\hfill$\square$

\textbf{Remark 2.9. The completion of the proof of Theorem 1.3.} It follows from Lemma 2.8 that, without loss of generality, we can assume that $G(A_A)=0$, i.e., $A$ is a right non-singular ring. In this case, Theorem 1.3 follows from Theorem 2.7.

\end{document}